\documentclass[11pt]{amsart}
\usepackage{graphicx}
\usepackage{amsmath}

\newtheorem{thm}{Theorem}[section]

\newtheorem{lem}[thm]{Lemma}

\theoremstyle{definition}

\theoremstyle{remark}
\newtheorem{rem}[thm]{Remark}
\numberwithin{equation}{section}
\begin{document}

\title[Effective bandwidth problem]
{The effective bandwidth problem revisited}%
\author{\bf Vyacheslav M. Abramov}%
\address{School of Mathematical Sciences, Monash University, Building 28M,
Clayton Campus, Clayton, VIC 3800, Australia}%
\email{vyacheslav.abramov@sci.monash.edu.au}%

%\thanks{}%

%\date{}%
%\dedicatory{}%
%\commby{}%
% ----------------------------------------------------------------
\begin{abstract}
The paper studies a single-server queueing system with autonomous
service and $\ell$ priority classes. Arrival and departure
processes are governed by marked point processes. There are $\ell$
buffers corresponding to priority classes, and upon arrival a unit
of the $k$th priority class occupies a place in the $k$th buffer.
Let $N^{(k)}$, $k=1,2,\ldots,\ell$ denote the quota for the total
$k$th buffer content. The values $N^{(k)}$ are assumed to be
large, and queueing systems both with finite and infinite buffers
are studied. In the case of a system with finite buffers, the
values $N^{(k)}$ characterize buffer capacities.
 The paper discusses a
circle of problems related to optimization of performance measures
associated with overflowing the quota of buffer contents in
particular buffers models. Our approach to this problem is new,
and the presentation of our results is simple and clear for real
applications.
\end{abstract}
\maketitle

%\tableofcontents
\newpage
% ----------------------------------------------------------------
\section{Introduction}\label{sec1}

\subsection{Approach}\label{sec1.1}

During the last two decades there has been an increasing interest
in the effective bandwidth problem for queueing systems with
priorities. There are different classes of messages (units)
arriving in telecommunication systems, and all of them are
characterized by their quality of service requirements. In order
to provide these quality of service guarantees and to allocate
necessary network resources, different priority classes
characterizing units arriving to that network are used.

There are a large number of papers related to this subject. A
detailed review of the related literature (up to publication time)
can be found in Berger and Whitt \cite{Berger and Whitt 1998a}
(for further discussions see also \cite{Berger and Whitt 1998b}).
For other relevant contributions to this subject see
\cite{Bestimas et al 1998}, \cite{Courcoubetis Siris Stamoulis
1999}, \cite{Elwalid and Mitra 1999}, \cite{Evans and Everitt
1999}, \cite{Kumaran et al 2000}, \cite{Lee et al 2005},
\cite{Wischik 1999}.

These papers all discuss approximations and suggest algorithms for
optimal solutions for the allocation of resources or effective
bandwidth problems. Most of these papers use large deviation
techniques. For example, Elwalid and Mitra \cite{Elwalid and Mitra
1995}, \cite{Elwalid and Mitra 1999} use Chernoff's inequality to
approximate loss probabilities in finite buffer systems with large
buffers. Berger and Whitt \cite{Berger and Whitt 1998a} also use
exponential asymptotics \cite{Berger and Whitt 1998b}, \cite{Whitt
1993} for the workload high level crossing of the $i$th class
priority unit. Other papers (e.g. \cite{Botvich and Duffield
1995}, \cite{Kelly 1996}, \cite{Paschalidis 1996}, \cite{Wischik
1999}, \cite{Wischik 2001}) also apply one or other techniques of
the large deviation principle. Many of the aforementioned papers
are aimed at solving concrete analytic problems, and their results
are based on an analysis of analytic transformations (such as
Laplace-Stieltjes or the $z$-transform) and their approximations.
Many of these results are then applied to $M/G/1$-oriented
queueing models or to models with more general arrival processes
having a Markov structure.

The approach of the present paper substantially differs from these
previous ones. The main focus of this paper is the solution of
bandwidth problems for $GI/M/1$-related priority systems. To the
knowledge of the author, such priority systems are not presented
in the literature where the overwhelming majority of priority
queueing systems studied are of $M/GI/1$ type. The innovations of
the present paper are as follows.
\smallskip

1. We consider models of queues with an \textit{autonomous service
mechanism} (see e.g. \cite{Borovkov 1976}, \cite{Borovkov 1984} as
well as Section \ref{sec1.3} of this paper). The main results of
our analysis are based on stochastic equations, and our models are
studied under a rather general setting and can be applied to a
broad class of real telecommunication systems. The obtained
stochastic equations are then used for analysis of particular
systems with exponentially distributed service times, which are a
subclass of queues with an autonomous service mechanism. Note that
martingale techniques for priority queueing systems (different
from the systems considered here) have been developed by Kella
\cite{Kella 1993}. However, the approach of \cite{Kella 1993}
differs from the
 present one. Specifically, \cite{Kella 1993} studies fluid
networks of parallel queues with dependent L\'evy inputs. It shows
that the special construction given in the paper can be applied to
the analysis of workload processes in $M/G/1$ queues with a
preemptive resume discipline. The paper of Kella \cite{Kella 1993}
is based on an extension of the earlier results of Kella and Whitt
\cite{Kella and Whitt 1992}.
 In contrast, our approach is based on a direct
construction of queues with autonomous service mechanisms, and can
be applied both to $M/G/1$ and $GI/M/1$ oriented priority queueing
systems. ($GI/M/1$ queues are precisely described in the paper.)

\smallskip

2. The buffer content process is described by the so-called
\textit{buffer type stochastic equation}. The buffer type
stochastic differential equation is a special stochastic
differential equation with discontinuous right-hand side (see
Filippov \cite{Filippov 1988}) and has already been used by
Elwalid and Mitra \cite{Elwalid and Mitra 1995}, \cite{Elwalid and
Mitra 1999} to study the model with two priority classes. However
\cite{Elwalid and Mitra 1995} and \cite{Elwalid and Mitra 1999}
used the explicit forms of this equation related to low and high
priority units. Analysis of these explicit equations is a hard
problem. In contrast, our buffer type equation are represented in
an (equivalent) integral form, and we discover a very simple
representation for cumulative buffer contents, see Lemma
\ref{lem1} and Theorem \ref{thm1}. According to this
representation, the system of equations for cumulative buffer
content processes is the usual system of stochastic equations
describing standard queue-length processes with an autonomous
service mechanism. This finding essentially simplify the analysis,
algorithms of solution and finally gives very simple approximation
of the explicit solution. For example, it enables us to study the
system with an arbitrary number of priorities.

\smallskip

3. Some papers (e.g. \cite{Elwalid and Mitra 1995}, \cite{Elwalid
and Mitra 1999}) assume that buffers have large capacities and
discuss the probabilities of buffer overflow. They use general
estimates given by large deviation theory, and particularly, by
Chernoff's inequality. Being well-motivated theoretically, these
estimates do not properly solve real practical problems. There is
an example in \cite{Choudhury Lucantoni and Whitt 1996} showing
that inequalities based on exponential bounds can give unrealistic
results.

We offer a \textit{unified approach} to systems with finite and
infinite buffers. Large parameters $N^{(k)}$, $k=1,2,\ldots,\ell$,
that are used in the sequel, are referred to as \textit{quota} for
buffer content and are related to finite and infinite buffers
systems. In the case of finite buffers models with recurrent input
and exponentially distributed service times of batches, we develop
the known asymptotic results on losses in $GI/M/1/n$ queues as
$n\to\infty$ \cite{Abramov 2002} to the case of $GI/M^{Y=C}/1/N_k$
queues ($k=1,2,\ldots,\ell$) with large buffers $N_k$ (the second
position $M^{Y=C}$ of the notation $GI/M^{Y=C}/1/N_k$ means that
the service time of units is exponentially distributed, and batch
size is equal to $C$) and then adapt the obtained asymptotic
result to estimate the loss probability in systems with large
finite buffers. The asymptotic representation of this paper, that
is used for the loss probability in $GI/M^{Y=C}/1/N_k$ queues and
then for the probability of buffer overflow, is preferable to
general type estimates such as Cramer or Chernoff inequalities.
The asymptotic results for the loss probability in
$GI/M^{Y=C}/1/N_k$ queues are expressed via the roots of the
appropriate functional equations (see Sections \ref{sec7} and
\ref{sec8}). They are also useful in studying the behaviour of
losses in the case of heavy load conditions. Cramer and Chernoff
inequalities are rougher, but their advantage is that they are
explicit. However there are exact estimates in the form of
explicit inequalities in the literature for the stationary
probabilities of $GI/M/1/N$ large buffer queueing systems as well
(see, \cite{Choi and Kim (2000)}), and they can be easily adapted
to the loss probabilities of the standard $GI/M/1/n$ queueing
system (see the discussion section in \cite{Abramov 2002}) with
application to models such as the $GI/M^{Y=C}/1/N_k$ queues
considered in the paper. We however are not going so far.

The results of this paper can be also applied to $M/G^{Y=C}/1/N$
oriented large buffers models. However, in this case a special
asymptotic analysis similar to that given in Abramov \cite{Abramov
2004 SIAM}  is necessary. This asymptotic analysis is routine and
not provided in the paper.

\subsection{Convention on the notation}\label{sec1.2}
For any increasing random sequence of points $t_1$, $t_2$,\ldots,
the associated point process
$Z(t)=\sum_{n=1}^\infty\mathbf{I}\{t_n\leq t\}$ is always denoted
by a capital Latin letter. If $\zeta_1$, $\zeta_2$, \ldots is a
sequence of marks, then the associated marked point process
$\mathcal{Z}(t)=\sum_{j=1}^{Z(t)}\zeta_j$ is always denoted by
calligraphic letters. All processes considered in the paper are
assumed to be right-continuous having left-limits and starting at
zero. Exceptions from this rule are especially mentioned in the
text (e.g. Remark \ref{rem1}). For an arbitrary point process
$\mathcal{Z}(t)$, its jump in point $t$ is denoted $\triangle
\mathcal{Z}(t)=\mathcal{Z}(t)-\mathcal{Z}(t-)$, where
$\mathcal{Z}(t-)$ is the left-limit of the process in point $t$.
For arrival processes we use letters $A$ and $\mathcal{A}$ with
sub- or super-script (the notation is given in Section
\ref{sec1.3}), and for departure process we use letters $D$ and
$\mathcal{D}$. The buffer processes describing the buffer contents
will be denoted by calligraphic letter $\mathcal{Q}$ with sub- or
super-script (the notation is in Section \ref{sec1.3}).  All
processes of this paper are assumed to be given on a common
filtered probability space $\{\Omega, \mathcal{F},
\mathbf{F}=(\mathcal{F}_t)_{t\geq0}, \mathbb{P}\}$.

\subsection{Description of the system}\label{sec1.3}
The paper is concerned with priority queueing system having $\ell$
buffers. Units arrive at the $k$th buffer at random time instants
$t_1^{(k)}=\tau_1^{(k)}$, $t_2^{(k)}=\tau_1^{(k)}+\tau_2^{(k)}$,
\ldots, and the $n$th unit arriving at the $k$th buffer has a
positive integer random length $\vartheta_n^{(k)}$. (In
telecommunication systems length can represent required memory for
the message.) Denote
$A^{(k)}(t)=\sum_{n=1}^\infty\mathbf{I}\left\{t_n^{(k)}\leq
t\right\}$, where $\mathbf{I}\{\cdot\}$ denotes an indicator of
the event, and
$\mathcal{A}^{(k)}(t)=\sum_{j=1}^{A^{(k)}(t)}\vartheta_j^{(k)}$.

The departure process $\mathcal{D}(t)$ is assumed to be a point
process with constant positive integer jumps $C$. Let $\chi_1$,
$\chi_2$,\ldots denote times between departures, and let
$x_n=\sum_{i=1}^n\chi_i$ denote the $n$th departure moment. Then
$\mathcal{D}(t)=C\sum_{n=1}^\infty \mathbf{I}\{x_n\leq t\}$. The
constant $C$ is called \emph{depletion rate}.

The buffers are numbered 1,2,\ldots,$\ell$, and the buffer with
lower order number has higher priority. Assume that the buffers
are infinite. Then the equation for the first buffer content
(\textit{highest} priority buffer) is
\begin{equation}
\label{1.1}\mathcal{Q}^{(1)}(t)=\max\left\{0,
\mathcal{Q}^{(1)}(t-)+\triangle \mathcal{A}^{(1)}(t)-\triangle
\mathcal{D}(t)\right\}.
\end{equation}
According to \eqref{1.1}, the buffer content
$\mathcal{Q}^{(1)}(t)$ is governed by the processes
$\mathcal{A}^{(1)}(t)$ and $\mathcal{D}(t)$ and is referred to as
a queueing process with an \textit{autonomous service mechanism}.
Queues with autonomous service mechanism were introduced and
originally studied by Borovkov \cite{Borovkov 1976},
\cite{Borovkov 1984}. For different applications see \cite{Abramov
2000}, \cite{Abramov 2004}, \cite{Abramov 2005}, \cite{Abramov
2008}, \cite{Fricker 1986}, \cite{Fricker 1987} and \cite{Gelenbe
and Iasnogorodski 1979}. The term $\triangle \mathcal{A}^{(1)}(t)$
is called the \emph{arrival jump} at time $t$, and the term
$\triangle \mathcal{D}(t)$ is called the \emph{possible departure
jump} at time $t$. The prefix \emph{possible} underlines the fact
that departures can occur only if the system is not empty. For
further simplifications, throughout the paper we assume that
arrival and departure processes are \textit{disjoint}, i.e. the
probability of simultaneous arrival and departure is 0.

If $t$ is a jump point of the process $\mathcal{D}(t)$, then the
\emph{real
 departure jump} at time $t$ is $\min\{\mathcal{Q}^{(1)}(t-)$, $C\}$.
 Thus, if $\mathcal{Q}^{(1)}(t-)$=0, then there is no departure jump.

$\mathcal{Q}^{(2)}(t)$ is the second buffer content, the priority
of which is lower than that of the first buffer.
$\mathcal{Q}^{(2)}(t)$ satisfies the equation:
\begin{multline}\label{1.2}
\mathcal{Q}^{(2)}(t) = \max\Big\{0, \mathcal{Q}^{(2)}(t-)+\triangle \mathcal{A}^{(2)}(t)\\
-\Big[\triangle
\mathcal{D}(t)-\mathcal{Q}^{(1)}(t-)\Big]\mathbf{I}\{\mathcal{Q}^{(1)}(t)=0\}\Big\}.
\end{multline}
Despite the fact that equation \eqref{1.2} has a more complicated
form than equation \eqref{1.1}, both of these equations are of the
same type. The term $\triangle \mathcal{A}^{(2)}(t)$ is an arrival
jump at time $t$. The structure of the departure jump is more
difficult. For simplicity we discuss the case $\ell=2$ below.  If
$t$ is at a departure jump and $\{\mathcal{Q}^{(1)}(t)>0\}$, then
$\{\mathcal{Q}^{(1)}(t-)> C\}$, and the jump is related to the
first buffer only. Otherwise, if $\{\mathcal{Q}^{(1)}(t)=0\}$,
then the following two cases are possible:
\begin{equation*}
\left\{0<\mathcal{Q}^{(1)}(t-)\leq C\right\},\leqno(\mbox{i})
\end{equation*}
\begin{equation*}
\left\{\mathcal{Q}^{(1)}(t-)=0\right\}.\leqno(\mbox{ii})
\end{equation*}
In case (i) departures occur from the first buffer, the first
buffer is completely emptied, and if the second buffer is not
empty, then in the case $\{\mathcal{Q}^{(1)}(t-)< C\}$ departures
occur also from the second buffer. In case (ii) departures occur
merely from the second buffer, provided that this buffer is not
empty. Thus the real departure jump in this case is
\begin{equation*}
\min\left\{\mathcal{Q}^{(1)}(t-)+\mathcal{Q}^{(2)}(t-), C\right\}.
\end{equation*}

Equation \eqref{1.2} is easily extended to the $k$th buffer
content for any $k$ =1,2,\ldots,$\ell$. Indeed, denoting
\begin{eqnarray}
\mathcal{Q}_k(t)&=&\mathcal{Q}^{(1)}(t)+\mathcal{Q}^{(2)}(t)+\ldots
+\mathcal{Q}^{(k)}(t),\label{1.2+}\\
\mathcal{A}_k(t)&=&\mathcal{A}^{(1)}(t)+\mathcal{A}^{(2)}(t)+\ldots+\mathcal{A}^{(k)}(t),\label{1.2++}
\end{eqnarray}
we have the following equation ($k=1,2,\ldots,\ell-1$):
\begin{multline}\label{1.3}
\mathcal{Q}^{(k+1)}(t) = \max\Big\{0, \mathcal{Q}^{(k+1)}(t-)+\triangle \mathcal{A}^{(k+1)}(t)\\
-\Big[\triangle
\mathcal{D}(t)-\mathcal{Q}_k(t-)\Big]\mathbf{I}\{\mathcal{Q}_k(t)=0\}\Big\}.
\end{multline}

The extension of \eqref{1.2} given by \eqref{1.3} is quite clear.
The term $\triangle \mathcal{A}^{(k+1)}(t)$ is an arrival jump  at
time $t$ (if any) to the buffer content $\mathcal{Q}^{(k+1)}(t-)$.
The other term of \eqref{1.3}
\begin{equation*}
\Big[\triangle
\mathcal{D}(t)-\mathcal{Q}_k(t-)\Big]\mathbf{I}\{\mathcal{Q}_k(t)=0\}
\end{equation*}
is also similar to the corresponding term of \eqref{1.2}. If $t$
is a jump point, then the meaning of $\mathcal{Q}_k(t-)$ is the
total content of all  buffers, the priority of which is greater
than the priority of the given $k+1$st buffer before the jump at
point $t$, and $\{\mathcal{Q}_k(t)=0\}$ is the event, that all
buffers, the priority of which is greater than the priority of the
given $k+1$st buffer, are empty after the jump at time $t$.

In the sequel the process $\mathcal{Q}_k(t)$ is called the $k$th
\textit{cumulative} buffer content.

\subsection{Formulation of the problems}\label{sec1.4}
 The paper is concerned with the following problems.
Let $N^{(1)}$, $N^{(2)}$, \ldots $N^{(\ell)}$ be large positive
integer values. Assuming that appropriate limits in probability
exist, denote
\begin{eqnarray}\label{1.4}
J^{(k)}=\mathbb{P}^{\_}\lim_{t\to\infty}\frac{1}{{A}_\ell(t)}
\sum_{j=1}^{{A}^{(k)}(t)}\mathbf{I}\left\{\mathcal{Q}^{(k)}(t_j^{(k)})>
N^{(k)}\right\},\\
k=1,2,\ldots,\ell.\nonumber
\end{eqnarray}
$A_\ell(t)$ is the total number of arrivals until time $t$. Then
$J^{(k)}$ is the fraction of arrival instants when the length
$N^{(k)}$ of the $k$th buffer is exceeded. Let $\alpha^{(1)}$,
$\alpha^{(2)}$,\ldots, $\alpha^{(\ell)}$ be real positive numbers,
denoting cost rates, and
\begin{equation}\label{1.4+}
J=\alpha^{(1)}J^{(1)}+\alpha^{(2)}J^{(2)}+\ldots+\alpha^{(\ell)}
J^{(\ell)}.
\end{equation}

Typical questions arising here are the following.

1. Assume that the parameters $N^{(1)}$, $N^{(2)}$, \ldots,
$N^{(\ell)}$ are given, but the depletion rate $C$ can be
controlled.
 Under what value of the depletion rate $C$ we have $J\leq \varepsilon$,
 where $\varepsilon$ is a given positive small value? This
 question can be formally written as follows: minimize $C$ subject
 to $J\leq \varepsilon$.
 \smallskip

2. Assume that $C$ is given, but $N^{(1)}$, $N^{(2)}$, \ldots,
$N^{(\ell)}$ are control variables. Assume additionally that with
given $\beta^{(2)}$, $\beta^{(3)}$,\ldots,$\beta^{(\ell)}$ the
values $N^{(1)}$, $N^{(2)}$, \ldots, $N^{(\ell)}$ must satisfy the
condition: $N^{(1)}=\lfloor\beta^{(2)}
N^{(2)}\rfloor$=$\lfloor\beta^{(3)} N^{(3)}\rfloor$ = \ldots =
$\lfloor\beta^{(\ell)} N^{(\ell)}\rfloor$, where $\lfloor
\cdot\rfloor$ is the notation for the integer part of number. The
problem is to minimize $N^{(1)}$ subject to $J\leq \varepsilon$.

\smallskip

\begin{rem}\label{rem1}
\eqref{1.4} applies to the finite and infinite buffers systems. To
finite buffers model we prescribe that a complete  \textit{
arrival group is rejected} when upon arrival the buffer overflows.
In the case of the system with infinite buffers,
$\mathcal{Q}^{(k)}(t)$, $k=1,2,\ldots,\ell$ all are assumed to be
right continuous having left limits. In the case of finite buffers
model, $\mathcal{Q}^{(k)}(t)$, $k=1,2,\ldots,\ell$ are not longer
right-continuous. For example, if $N^{(1)}$ is the capacity of the
first buffer, and at moment $t_j^{(1)}$ the buffer overflows, then
we admit that $\mathcal{Q}^{(1)}\left(t_j^{(1)}\right)$ is greater
than $N^{(1)}$ in $t_j^{(1)}$. However in the neighborhood of this
point $\mathcal{Q}^{(1)}\left(t_j^{(1)}\right)\leq N^{(1)}$. Then
the left and right limits of $\mathcal{Q}^{(1)}(t)$ in point
$t_j^{(1)}$ are not greater than $N^{(1)}$, both these limits
(with probability 1) are equal and $t_j^{(1)}$ is an isolated
point.
\end{rem}

\subsection{Brief description of the mathematical ideas,
methodology and contribution of the paper}\label{Ideas} In this
section we describe the mathematical ideas of this paper, as well
as the methodology and overall contribution.

We start from the description of the buffer content process. For
the highest priority buffer the equation for the buffer content is
very simple. It is described by equation \eqref{1.1}. The
equations for the lower priority buffer contents are relatively
more complicated and described by equation \eqref{1.2}. However,
the equations for the cumulated buffer contents are simple and
described by a difference recurrence equation similar to
\eqref{1.1}
\begin{equation}
\label{I1}\mathcal{Q}_k(t)=\max\{0, \mathcal{Q}_k(t-)+\triangle
\mathcal{A}_k(t)-\triangle \mathcal{D}(t)\}.
\end{equation}
Another form for \eqref{I1} is a stochastic equation
\begin{equation}
\label{I3}\mathcal{Q}_k(t)=\mathcal{A}_k(t)-\sum_{j=1}^C\int_0^t
\mathbf{I}\{\mathcal{Q}_k(u-)\geq j\} \mathrm{d}D(u).
\end{equation}
(In all these two equations $k=1,2,\ldots,\ell$.) The stochastic
equation \eqref{I3} can be rewritten
\begin{equation*}
\mathcal{Q}_k(t)=\mathcal{A}_k(t)-\mathcal{D}(t)+\sum_{j=1}^C\int_0^t\mathbf{I}\{\mathcal{Q}_k(u-)\leq
j-1\}\mathrm{d}D(u),
\end{equation*}
with subsequent reduction to a Skorokhod problem (see
\cite{Anulova and Liptser 1990}, \cite{Skorokhod 1961},
\cite{Tanaka 1979}). (In the case $C=1$ such reduction was
provided in \cite{Kogan and Liptser 1993}. For its further
application see also \cite{Abramov 2000} and \cite{Abramov 2004}.)

It is shown then that representations similar to \eqref{I1} and
\eqref{I3} remain valid for finite buffer models. Thus, in all
cases the problem reduces to analyzing queueing systems with an
autonomous service mechanism.

We use these results for analysis of particular queueing buffer
models with priorities. A system with exponentially distributed
service times is a special case of a system with an autonomous
service mechanism.  (A special construction of models with finite
and infinite buffers is explained later in Section \ref{sec7} of
the paper.)

In the case of finite buffer models with renewal input and
exponentially distributed service times we adapt recent results on
asymptotic analysis \cite{Abramov 2002}. As in papers
\cite{Abramov 2002}, \cite{Abramov 2004 SIAM} and \cite{Abramov
2006}, the analysis is based on reducing the loss probability to a
convolution type recurrence relation:
\begin{equation*}
f_n = \sum_{j=0}^n f_{n-j+1}\pi_j \ \ (f_0>0),
\end{equation*}
where $\pi_0>0$, $\pi_j\geq0$ for all $j=1,2,\ldots$, and
$\sum_{j=0}^\infty\pi_j=1$, and applying asymptotic analysis
similar to that of the book of Tak\'acs \cite{Takacs 1967},
p.22-23. Consequently, we provide heavy traffic analysis of these
models based on asymptotic expansions of the results obtained
under ``usual" conditions. The loss probability for the large
finite $k$th cumulative buffer is then not greater than the sum of
the loss probabilities in the associated $GI/M^{Y=C}/1/N_i$
queues, $i=1,2,\ldots,k$. However, for large values $N_i$ this sum
is very small, with the order of this sum being the same as the
order of one (maximum) term obtained by asymptotic analysis, and
an estimate obtained seems to be better than that estimate
obtained by rough methods of large deviation principle and
Chernoff's inequality.

\smallskip

Thus, the main mathematical contribution is a general theory of
priority buffer models with application to particular priority
queueing systems with recurrent input and large buffers.

\subsection{Organization of the paper}
In Section \ref{sec3}, Lemma \ref{lem1} states that the $k$th
cumulative buffer content has the representation \eqref{2.0}.  The
intuitive sense of Lemma \ref{lem1} is that the $k$th cumulative
buffer content for the system with infinite buffers is described
by the same equation as the queue-length process in the queueing
system with autonomous service mechanism, the arrival process of
which is $\mathcal{A}_k(t)$ and the departure process
$\mathcal{D}(t)$. We further prove a stability theorem. The main
condition for stability is \eqref{3.5}, the proof being based on
reduction to the Skorokhod reflection principle and results of
Borovkov \cite{Borovkov 1976}, \cite{Borovkov 1984}. In section
\ref{sec2add} the model with finite buffers is considered. It is
shown that the equations for cumulative buffer contents in this
case are similar to the case of a model with infinite buffers. In
Section \ref{sec4} we derive the formula for $J^{(k)}$,
$k=1,2,\ldots,\ell$, using the level-crossing method based on
representation \eqref{4.3}. In Sections  \ref{sec7} and \ref{sec8}
special models of queueing systems are studied. The results of
these sections are illustrative, and we do not discuss general
buffer models with batch arrival such as $GI^X/M^{Y=C}/1$ queues,
although the asymptotic geometrical bounds for stationary
probability to reach high level $N$ in $GI^X/M^Y/1$ queues is
known (see \cite{Economou and Fakinos 2003}).  All models
considered here are particular cases of the general models
discussed in Sections \ref{sec3} and \ref{sec2add}: these models
are with independent identically distributed interarrival times.
The results of Section \ref{sec8} are based on an extension of
recent results \cite{Abramov 2002}. As in \cite{Abramov 2002} the
asymptotic analysis is based on reduction to appropriate
representation helping us to use then the Tak\'{a}cs theorem on
asymptotic behavior of the convolution type recurrence relation
\cite{Takacs 1967}, p. 22-23.  In Section \ref{sec8.1} the
asymptotic behaviour of losses are studied under ``usual"
conditions, while in Section \ref{sec8.2} the analysis of losses
is done under heavy load conditions. In Section \ref{sec9}
approximation of the initial problem stated in Section
\ref{sec1.4} by another related problem is suggested. In Section
\ref{sec10} algorithms for numerical solution of the problems of
Section \ref{sec9} are proposed. There are concluding remarks in
Section \ref{sec12}.

\section{The stability theorem for the infinite buffers system}\label{sec3}
The representation for the buffer content of infinite buffers
systems given by \eqref{1.1}, \eqref{1.2} and \eqref{1.3} is
difficult to analyze. However, for the cumulative buffer contents
of infinite buffers systems the representation is simple.

\begin{lem} \label{lem1} For all $k=1,2,\ldots,\ell$ the following
equation for the $k$th cumulative buffer content
$\mathcal{Q}_k(t)$ holds:
\begin{equation}
\label{2.0}\mathcal{Q}_k(t)=\mathcal{A}_k(t)-\sum_{j=1}^C\int_0^t
\mathbf{I}\{\mathcal{Q}_k(u-)\geq j\} \mathrm{d}D(u),
\end{equation}
where $D(t)=\frac{\mathcal{D}(t)}{C}$.
\end{lem}

The proof of this lemma is given in Appendix A.

The statement of Lemma \ref{lem1} has a simple intuitive
explanation. For example, in the case $\ell=2$ we have two classes
of units, and clearly the cumulative buffer content process
$\mathcal{Q}_2(t)=\mathcal{Q}^{(1)}(t)+\mathcal{Q}^{(2)}(t)$
contains two unit classes together, and therefore must behave as a
usual (i.e. without priorities) queue-length process with an
autonomous service mechanism, the arrival process of which is
$\mathcal{A}_2(t)=\mathcal{A}^{(1)}(t)+\mathcal{A}^{(2)}(t)$, and
the departure process is $\mathcal{D}(t)$. This intuitive
explanation is easily extended to the case of arbitrary
$k=1,2,\ldots,\ell$ number of classes.

The right-hand side of this equation contains the  sum
$$
\sum_{j=1}^C\int_0^t \mathbf{I}\{\mathcal{Q}_k(u-)\geq j\}
\mathrm{d}D(u).
$$
Nevertheless, the problem can be reduced to the Skorokhod
reflection principle.

Denote $S_k(t)=\mathcal{A}_k(t)-\mathcal{D}(t)$,
$k=1,2,\ldots,\ell$. Then,
\begin{equation}\label{Sk1}
\mathcal{Q}_k(t)=S_k(t)+\sum_{j=1}^C\int_0^t\mathbf{I}\{\mathcal{Q}_k(u-)\leq
j-1\}\mathrm{d}D(u).
\end{equation}
Equation \eqref{Sk1} implies that $\mathcal{Q}_k(t)$ is the normal
reflection of the process $S_k(t)$ ($S_k(0)=0$) at zero. More
accurately, $\mathcal{Q}_k(t)$ is the nonnegative solution of the
Skorokhod problem of the normal reflection of the process $S_k(t)$
at zero (see Skorokhod \cite{Skorokhod 1961} as well as Tanaka
\cite{Tanaka 1979} and Anulova and Liptser \cite{Anulova and
Liptser 1990}, Ramanan \cite{Ramanan 2006}). This is because the
function
\begin{equation*}
\phi_k(t)=\sum_{j=1}^C\int_0^t\mathbf{I}\{\mathcal{Q}_k(u-)\leq
j-1\}\mathrm{d}D(u)
\end{equation*}
satisfies the following two properties:

(a) $\int_0^th[\mathcal{Q}_k(u)]\mathrm{d}\phi_k(t)=0$ for any
continuous nonnegative function $h(x)$ with $h(0)=0$;
\smallskip

(b) the function
$\int_0^t[Y(u)-\mathcal{Q}_k(u)]\mathrm{d}\phi_k(u)$ is not
decreasing for any nonnegative right-continuous function $Y(u)$
having the left limits.

\smallskip
Let us show (a). We have
\begin{equation*}
\begin{aligned}
\int_0^th[\mathcal{Q}_k(u)]\mathrm{d}\phi_k(u)&=\int_0^t
h[\triangle\mathcal{Q}_k(u)]\mathrm{d}\phi_k(u)\\
&=\int_0^t
h[\triangle\mathcal{Q}_k(u)]\mathrm{d}\left(\sum_{j=1}^C\int_0^u\mathbf{I}\{\mathcal{Q}_k(v-)\leq
j-1\}\mathrm{d}D(v)\right).
\end{aligned}
\end{equation*}
Let $u_i$ denote the points of jump of the process
$\mathcal{Q}_k(u)$ in the interval [$0,t$]. For the last integral
we have the following representation:
\begin{equation*}
\begin{aligned}
&\int_0^t
h[\triangle\mathcal{Q}_k(u)]\mathrm{d}\left(\sum_{j=1}^C\int_0^u\mathbf{I}\{\mathcal{Q}_k(v-)\leq
j-1\}\mathrm{d}D(v)\right)\\
&=\sum_{0\leq u_i\leq t}\sum_{j=1}^C
h[\triangle\mathcal{Q}_k(u_i)]\triangle
\mathcal{A}_k(u_i)\mathbf{I}\{\mathcal{Q}_k(u_i-)\leq
j-1\}\triangle D(u_i).
\end{aligned}
\end{equation*}
The last sum is a finite sum: the number of points $u_i$ is finite
in any finite interval [$0,t$] with probability 1. Any value of
jump $\triangle\mathcal{Q}_k(u_i)$ is bounded with probability 1,
and the nonnegative continuous function
$h[\triangle\mathcal{Q}_k(u)]$, satisfying the property $h(0)=0$
is therefore bounded for all $0\leq u\leq t$. In addition, taking
into account that the jumps of the processes $\mathcal{A}_k(u)$
and $D(u)$ are disjoint, i.e. either
$\triangle\mathcal{A}_k(u_i)=0$ or $\triangle{D}(u_i)=0$ with
probability 1, we arrive at the conclusion that
$\int_0^th[\mathcal{Q}_k(u)]\mathrm{d}\phi_k(u)=0$.
 (a) follows.

 (b) is
implied by (a).

It follows from the Skorokhod reflection principle that the
function $\phi_k(t)$ has the following representation:
\begin{equation*}
\phi_k(t)=-\inf_{u\leq t} S_k(u).
\end{equation*}

Therefore $\mathcal{Q}_k(t)$ has the following representation
\begin{equation}\label{3.4}
\begin{aligned}
&\mathcal{Q}_k(t)=S_k(t)-\inf_{u\leq t}S_k(u),\\
&\ \ \ \ \ \ k=1,2,\ldots,\ell.
\end{aligned}
\end{equation}

Equation \eqref{3.4} is well-known in queueing theory. Following
Borovkov \cite{Borovkov 1976}, we have the following statement of
the stability.

\begin{thm}\label{thm1} Assume
\begin{equation}
\label{3.5} \mathbb{P}\left\{\lim_{t\to\infty}
\frac{\mathcal{A}_\ell(t)-\mathcal{D}(t)}{t}=r<0\right\},
\end{equation}
and $\widetilde S_k(t)$, $k=1,2,\ldots,\ell$, are stationary point
processes, the increments of which coincide in distribution with
the corresponding increments of the processes
$\mathcal{A}_k(t)-\mathcal{D}(t)$, $k=1,2,\ldots,\ell$.

Then there exist  stationary processes $\mathcal{Q}^{(k)}(T)$,
$k=1,2,\ldots,\ell$, such that
\begin{equation}\label{3.6}
\mathcal{Q}^{(1)}(T){\buildrel 'd'\over =}\sup_{u\leq
T}\left[\widetilde S_1(T)-\widetilde S_1(u)\right],
\end{equation}
and
\begin{equation}\label{3.7}
\begin{aligned}
&\mathcal{Q}^{(k)}(T){\buildrel 'd'\over =}\sup_{u\leq
T}\left[\widetilde S_k(T)-\widetilde S_k(u)\right]-\sup_{u\leq
T}\left[\widetilde
S_{k-1}(T)-\widetilde S_{k-1}(u)\right]\\
&k=2,3,\ldots,\ell.
\end{aligned}
\end{equation}
\end{thm}
\begin{proof}
The proof is based on representation \eqref{3.4} and can be found
in Borovkov \cite{Borovkov 1976}. Specifically, it follows from
that proof that there are stationary processes $\mathcal{Q}_k(T)$,
$k=1,2,\ldots,\ell$ such that
\begin{equation}\label{3.8}
\mathcal{Q}_k(T){\buildrel 'd'\over =}\sup_{u\leq
T}\left[\widetilde S_k(T)-\widetilde S_k(u)\right].
\end{equation}
Therefore, keeping in mind that
$\mathcal{Q}^{(k)}(t)=\mathcal{Q}_k(t)-\mathcal{Q}_{k-1}(t)$,
$k=2,3,\ldots,\ell$ and $\mathcal{Q}^{(1)}(t)=\mathcal{Q}_1(t)$,
from \eqref{3.8} we have \eqref{3.6} and \eqref{3.7}.
\end{proof}

\section{The finite buffers model}\label{sec2add}
Equation \eqref{2.0} and other related equations for infinite
buffers content can be easily extended for the model with finite
buffers. It is assumed that if upon arrival of a batch the buffer
of a given class overflows, then the complete arrival batch is
rejected, see Remark \ref{rem1}.

For the analysis of the finite buffers case we introduce new
arrival processes $\overline{\mathcal{A}}^{(k)}(t)$, which are
derived from the initial processes $\mathcal{A}^{(k)}(t)$ as
follows. We set
\begin{equation}\label{4+.0}
\triangle \overline{\mathcal{A}}^{(k)}(t)=\triangle
\mathcal{A}^{(k)}(t)\mathbf{I}\left\{\mathcal{Q}^{(k)}(t)\leq
N^{(k)}\right\}.
\end{equation}
The arrival processes $\overline{\mathcal{A}}^{(k)}(t)$ take into
account only jumps of real buffer content process. Thus
$A(t)-\overline{A}(t)$ is the number of lost units during time
$t$, and $\mathcal{A}(t)-\overline{\mathcal{A}}(t)$ is their total
length during that time $t$.

Then the buffer content process $\mathcal{Q}^{(1)}(t)$ is defined
by the pair of equations
\begin{equation}
\label{4+.1}\mathcal{Q}^{(1)}(t)=\max\left\{0,
\mathcal{Q}^{(1)}(t-)+\triangle \mathcal{A}^{(1)}(t)-\triangle
\mathcal{D}(t)\right\},
\end{equation}
\begin{equation}
\label{4+.2}\mathcal{Q}^{(1)}(t+)=\max\left\{0,
\mathcal{Q}^{(1)}(t-)+\triangle
\overline{\mathcal{A}}^{(1)}(t)-\triangle \mathcal{D}(t)\right\}.
\end{equation}
Thus in the case $\triangle \mathcal{A}^{(1)}(t)=\triangle
\overline{\mathcal{A}}^{(1)}(t)$ the buffer contents
$\mathcal{Q}^{(1)}(t)$ and $\mathcal{Q}^{(1)}(t+)$ are equal and
there is no loss at time $t$. Otherwise, if $\triangle
\mathcal{A}^{(1)}(t)\neq\triangle
\overline{\mathcal{A}}^{(1)}(t)$, i.e. $\triangle
\overline{\mathcal{A}}^{(1)}(t)=0$ and $\triangle
\mathcal{A}^{(1)}(t)>0$, then there is a loss of a unit in time
$t$.

Next, similarly to \eqref{1.3} for $k=1,2,\ldots,\ell-1$ we have
another pair of equations:
\begin{multline}\label{4+.3}
\mathcal{Q}^{(k+1)}(t) = \max\Big\{0, \mathcal{Q}^{(k+1)}(t-)+\triangle \mathcal{A}^{(k+1)}(t)\\
-\Big[\triangle \mathcal{D}(t)-\mathcal{Q}_k(t-)-\triangle
\mathcal{A}_k(t)\Big]\mathbf{I}\{\mathcal{Q}_k(t)=0\}\Big\},
\end{multline}
\begin{multline}\label{4+.4}
\mathcal{Q}^{(k+1)}(t+) = \max\Big\{0, \mathcal{Q}^{(k+1)}(t-)+\triangle \overline{\mathcal{A}}^{(k+1)}(t)\\
-\Big[\triangle \mathcal{D}(t)-\mathcal{Q}_k(t-)-\triangle
\overline{\mathcal{A}}_k(t)\Big]\mathbf{I}\{\mathcal{Q}_k(t)=0\}\Big\}.
\end{multline}

Similarly to Lemma \ref{lem1}, for the finite buffers model we
have the following lemma.
\begin{lem}
\label{lem3} For all continuity points of the $k$th cumulative
buffer content process $\mathcal{Q}_k(t)$, $k=1,2,\ldots,\ell$, we
have:
\begin{equation}
\label{4+.5}\mathcal{Q}_k(t)=\overline{\mathcal{A}}_k(t)-\sum_{j=1}^C\int_0^t
\mathbf{I}\{\mathcal{Q}_k(u-)\geq j\} \mathrm{d}D(u).
\end{equation}
\end{lem}

\section{The formula for $J^{(k)}$}\label{sec4}

In this section we study the dynamics of the buffer lengths by
level-crossings analysis for the infinite buffers model. It is
assumed throughout that condition \eqref{3.5} for the stability is
fulfilled.

In addition to the stability condition assume:
\begin{eqnarray}\label{4.1}
\mathbb{P}\left\{\lim_{t\to\infty}\frac{{A}^{(k)}(t)}{t}=\lambda^{(k)}\right\}=1,\\
k=1,2,\ldots,\ell,\nonumber
\end{eqnarray}
and
\begin{equation}\label{4.2}
\mathbb{P}\left\{\lim_{t\to\infty}\frac{{D}(t)}{t}=\mu\right\}=1.
\end{equation}

Then according to \eqref{3.5} the sequences
$\frac{1}{n}\sum_{i=1}^{n}\vartheta_i^{(k)}$, $k=1,2$,\ldots,
$\ell$, as $n\to\infty$, also converges with probability 1.

Recall that $t_1^{(k)}=\tau_1^{(k)}$,
$t_2^{(k)}=\tau_1^{(k)}+\tau_2^{(k)}$, \ldots,
($k=1,2,\ldots,\ell$) denote the sequence of points (arrival
moments) of the process $\mathcal{A}^{(k)}(t)$, and $x_1=\chi_1$,
$x_2=\chi_1+\chi_2$, \ldots denote the sequence of points (the
moments of possible departure jumps) of $\mathcal{D}(t)$.

Then, for the number of up- and down-crossings for $m\geq1$ we
have the following equation:
\begin{equation}\label{4.3}
\begin{aligned}
&\sum_{i=1}^{{A}^{(k)}(t)}\mathbf{I}\left\{\mathcal{Q}^{(k)}\Big(t_{i}^{(k)}\Big)\geq
m, \ \mathcal{Q}^{(k)}\Big(t_{i}^{(k)}-\Big)<m\right\}\\
&=
\sum_{j=1}^{D(t)}\mathbf{I}\left\{m\leq\mathcal{Q}^{(k)}\Big(x_{j}-\Big)
\leq m-1+C\right\}\\& \ \ \ +\mathbf{I}
\left\{\mathcal{Q}^{(k)}(t)\geq m\right\}\\
&=
\sum_{l=1}^C\sum_{j=1}^{D(t)}\mathbf{I}\left\{\mathcal{Q}^{(k)}\Big(x_{j}-\Big)
=m-1+l\right\}\\& \ \ \ +\mathbf{I}
\left\{\mathcal{Q}^{(k)}(t)\geq m\right\},
\end{aligned}
\end{equation}
where $\mathcal{Q}^{(k)}(0)=0$. Equation \eqref{4.3} can be
explained as follows. The left-hand side of the equation is the
number of arrivals until time $t$, seeing before arrival the
buffer content less than $m$ and at the moment of arrival not
smaller than $m$. This constitutes the number of up-crossings of
the level $m$ until time $t$, i.e. the number of instants where
arrivals jump over the level $m-1$. The first term of the
right-hand side describes the number of departure moments when
immediately before departure the buffer content is between $m$ and
$m+C-1$. (Then after the departure the buffer content is between
$\max\{0, m-C\}$ and $m-1$, and this constitutes the number of
down-crossings of the level $m$). The difference between the
number of up-crossings and down-crossings of level $m$ can be
either 1 or 0, and the second term of the right-hand side
compensates for this difference.

Dividing the both sides of \eqref{4.3} by $t$, and letting $t$
increase unboundedly, we  obtain:
\begin{equation}
\label{4.4}
\begin{aligned}
&\lim_{t\to\infty}\frac{1}{t}\mathbb{E}
\sum_{i=1}^{{A}^{(k)}(t)}\mathbf{I}\left\{\mathcal{Q}^{(k)}\Big(t_{i}^{(k)}\Big)\geq
m, \ \mathcal{Q}^{(k)}\Big(t_{i}^{(k)}-\Big)<m\right\}\\ &=
\lim_{t\to\infty}\frac{1}{t}\mathbb{E}\sum_{l=1}^C\sum_{j=1}^{D(t)}\mathbf{I}\left\{\mathcal{Q}^{(k)}(x_{j}-)
=m-1+l\right\},
\end{aligned}
\end{equation}
and after elementary transformations (see Appendix B) we arrive at
\begin{equation}
\label{4.9}
\begin{aligned}
J^{(k)}
&=\frac{\lambda^{(k)}}{\lambda^{(1)}+\lambda^{(2)}+\ldots+\lambda^{(\ell)}}\cdot
\frac{1}{\lambda^{(k)}}\\
&\ \ \
\times\lim_{t\to\infty}\frac{1}{t}\mathbb{E}\int_0^t\sum_{l=1}^{C}\mathbf{I}
\left\{\mathcal{Q}^{(k)}(u-)\geq N^{(k)}+l\right\}\mathrm{d}{D}(u)\\
&=\frac{1}{\lambda_\ell}\lim_{t\to\infty}\frac{1}{t}\mathbb{E}\sum_{l=1}^{C}
\int_0^t\mathbf{I}
\left\{\mathcal{Q}^{(k)}(u-)\geq N^{(k)}+l\right\}\mathrm{d}{D}(u),\\
&\ \ \ \ \ \ k=1,2,\ldots,\ell,
\end{aligned}
\end{equation}
where $\lambda_\ell=\lambda^{(1)}+\lambda^{(2)}+\ldots
+\lambda^{(\ell)}$.

\section{The buffers content distribution of $GI/M^{Y=C}/1$
queues}\label{sec7}
\subsection{Main result}\label{MRGI-Infty}
We start this section with a representation for the buffer content
processes in the case where the arrival processes $A^{(k)}(t)$,
$k=1,2$,\ldots, $\ell$ all satisfy \eqref{4.1}, and the process
$D(t)$ is Poisson. Assume also that $\vartheta_n=1$ for all $n$.

We have:
\begin{equation}\label{7.1}
\begin{aligned}
&\lim_{t\to\infty}\frac{1}{t}\mathbb{E}\int_0^t\mathbf{I}\{\mathcal{Q}_k(u-)=m-1\}
\mathrm{d}A_k(u)\\&=\mu\lim_{t\to\infty}
\frac{1}{t}\sum_{l=1}^C\int_0^t\mathbb{P}\{\mathcal{Q}_k(u)=m-1+l\}
\mathrm{d}u,\\
& \ \ \ \ k=1,2,\ldots,\ell.
\end{aligned}
\end{equation}
However \eqref{7.1} does not permit us to obtain explicit results
for the stationary probabilities even in the case where the
processes $A^{(k)}(t)$, $k=1,2,\ldots,\ell$ all are renewal
processes. Moreover, in the case where all processes $A^{(k)}(t)$,
$k=1,2,\ldots,\ell$ are renewal, Lemma \ref{lem1} is no longer
useful in general, because stationary interarrival times to
cumulative buffers are dependent in general, and the corresponding
stationary arrival processes are not longer renewal.

Therefore we consider the following special case of the general
buffers model. Let $A(t)$ be a point process of arrivals
satisfying the condition
$\mathbb{P}\{\lim_{t\to\infty}\frac{A(t)}{t}$ = $\lambda\}=1$. Let
$\pi^{(1)}$, $\pi^{(2)}$,\ldots, $\pi^{(\ell)}$ be positive
probabilities, $\sum_{k=1}^\ell \pi^{(k)}=1$, where $\pi^{(k)}$ is
a probability that an arriving customer belongs to the class $k$.
Then the points processes $A^{(k)}(t)$, $k=1,2,\ldots,\ell$, are
all thinnings of the original process $A(t)$, and in the case
where $A(t)$ is a renewal process all the processes $A^{(k)}(t)$,
$k=1,2,\ldots,\ell$ are renewal processes as well with intensities
$\lambda^{(k)}=\lambda \pi^{(k)}$ correspondingly. Consequently,
the processes $A_k(t)$, $k=1,2,\ldots,\ell$, are renewal processes
with intensities $\lambda_k$ = $\lambda^{(1)}$ + $\lambda^{(2)}$
+\ldots+ $\lambda^{(k)}$, and one can apply the theory to the
$GI/M^{Y=C}/1$ and $GI/M^{Y=C}/1/N_k$ queues with large buffers
$N_k$.

By $GI/M^{Y=C}/1$ queue we mean a single-server queueing system
with recurrent input and exponentially distributed service time of
the constant size batch $C$. In the sequel we use the notation
$GI/M^{C}/1$ for these queueing systems. $GI/M^{C}/1$ queueing
systems are particular systems with an autonomous service
mechanism, and they are therefore described by buffer type
stochastic differential equations or by one of the above
equivalent forms of these equations. For these queueing systems
therefore Lemma \ref{lem1} remains true. Specifically, from this
lemma one can conclude that the cumulative buffer content
processes are described by the steady-state distributions of the
usual queue-length processes of the $GI/M^{C}/1$ queues. The
stability condition for these queues is
$\rho_\ell=\frac{\lambda_\ell}{\mu C}<1$.

Using a standard method, the limiting and stationary probabilities
for the cumulative buffer contents of $GI/M^{C}/1$ queues are
calculated as follows. Let $t_{k,j}$ denote the $j$th arrival
moment to one of the first $k$ buffers. Then for limiting and
stationary probability we have the following.
\begin{thm}\label{prop3}
For cumulative buffer contents $\mathcal{Q}_k(t)$,
$k=1,2,\ldots,\ell$.
\begin{equation}
\label{7.2}
P_{k,m}=\lim_{j\to\infty}\mathbb{P}\{\mathcal{Q}_k(t_{k,j}-)=m\}=\varsigma_k^m(1-\varsigma_k),
\end{equation}
where $\varsigma_k$ is the (unique) root of the functional
equation
\begin{equation}\label{7.3}
z=\widehat B_k(\mu-\mu z^C)
\end{equation}
in the interval (0,1), and $\widehat
B_k(s)=\int_0^\infty\mathrm{e}^{-sx}\mathrm{d}B_k(x)$ is the
Laplace-Stieltjes transform of the stationary distribution of
interarrival time $B_k(x)$ to the first $k$ buffers.
\end{thm}

\begin{proof}
The stationary probabilities of $GI^X/M^Y/1$ queues can be found
in Economou and Fakinos \cite{Economou and Fakinos 2003}, and the
statement for $GI/M^{C}/1$ queues can be deduced from their
result.\footnote{The following additional condition is missed in
the main statement of \cite{Economou and Fakinos 2003}: the common
divisor of possible values $X$ and $Y$ must be equal to 1.}
However, the direct proof of the result for the $GI/M^{C}/1$ queue
is much simpler than that reduction from the aforementioned
general result. Therefore below the direct proof of this theorem
is provided.

First of all notice, that according to \eqref{7.1} the state
probabilities immediately before arrival, $P_{k,m}$, are
\begin{equation*}
P_{k,m}=\lim_{t\to\infty}\frac{1}{\lambda_kt}\mathbb{E}\int_0^t\mathbf{I}\{\mathcal{Q}_k(u-)=m\}
\mathrm{d}A_k(u),
\end{equation*}
and $P_{k,m}=(1-z)z^m$ for some $z<1$. Let $f_m$ denote the number
of up- (down-) crossing of level $m$ during a busy period of
$GI/M^C/1$ queue (the number of cases where immediately before
arrival there are $m$ customers in the system). Then, by renewal
arguments $\mathbb{E}f_m=z^m$, and according to the total
expectation formula for any $m\geq 1$ we have the following
equation:
\begin{equation}\label{7.3.2}
z^m=\sum_{i=0}^\infty\int_0^\infty\mathrm{e}^{-\mu
x}\frac{\left(\mu xz^C \right)^i}{i!}~z^{m-1}\mathrm{d}B_k(x),
\end{equation}
where $B_k(x)$ is the probability distribution function of
interarrival time. Therefore, from \eqref{7.3.2} we obtain
\begin{equation*}\label{7.3.3}
\begin{aligned}
z&=\sum_{i=0}^\infty\int_0^\infty\mathrm{e}^{-\mu
x}\frac{\left(\mu
xz^C \right)^i}{i!}~\mathrm{d}B_k(x)\\
&=\widehat B_k(\mu-\mu z^C),
\end{aligned}
\end{equation*}
and the statement of Theorem \ref{prop3} follows. By standard
method (see e.g. \cite{Takacs 1962}, \cite{Gnedenko Kovalenko
1968}) one can prove that under the assumption $\rho_\ell<1$ there
exists a unique root of equation $z=\widehat B_k(\mu-\mu z^C)$ in
the interval (0,1).
\end{proof}

\subsection{Particular case} We consider an $M/M^C/1$ queueing system with infinite
buffers. This particular case is easily deduced from the statement
of Theorem \ref{prop3}. Specifically, in the case of Poisson
arrivals from \eqref{7.3} we obtain the equation:
\begin{equation*}
z=\frac{\lambda_k}{\lambda_k+\mu-\mu z^C}.
\end{equation*}
Then, the constant $\varsigma_k$ must be the solution of equation
\begin{equation}\label{7.3.4}
\frac{\lambda_k}{\mu}=\sum_{i=1}^Cz^i,
\end{equation}
belonging to the interval (0,1). A similar result can be also
found in \cite{Chao Pinedo and Shaw 1996} for nodes of a network,
the customers of which are served by random batches.

\section{Loss probabilities for cumulative buffers}\label{sec8}

In this section we discuss loss probabilities assuming that the
$k$th cumulative buffer content has large capacity $N_k$. We study
buffer loss probabilities under ``usual" and heavy load
conditions. By ``usual" conditions we mean the case when the load
parameter of the queueing system is fixed, while in the case of
heavy load conditions the sequence of load parameters, associated
with series of queueing systems, approaches 1.

\subsection{Loss probabilities under ``usual"
conditions}\label{sec8.1} We use the notation $GI/M^{C}/1/N_k$ for
the queueing systems with finite capacity $N_k$, similar to the
notation used for the queueing systems with infinite capacity in
the previous section. According to Lemma \ref{lem3} the cumulative
buffer contents in continuity points of the process
$\mathcal{Q}_k(t)$ behave as usual $GI/M^{C}/1/N_k$ queues.
However, the behavior of the number of losses, the main
characteristic of interest, is essentially different, that is the
losses in $GI/M^{C}/1/N_k$ queues are not equal to the losses in
the corresponding cumulative buffers $\mathcal{Q}_k(t)$.
Specifically, the losses in $GI/M^{C}/1/N_k$ queues occur only in
the case in which the buffer overflowed when the arriving customer
met all waiting places busy. The losses in the cumulative buffers
$\mathcal{Q}_k(t)$ can occur in many cases when one of specific
buffers, say $j$th buffer, $1\leq j\leq k$, has overflowed.

However, in some cases when the values $N_1<N_2<\ldots<N_\ell$ all
are large, a correspondence between $GI/M^{C}/1/N_k$ queues and
finite buffers models, may give useful asymptotic results.

Specifically, the loss probability of a customer arriving at one
of the first $k$ buffers is not greater than $p_1+p_2+\ldots+p_k$,
where $p_i$ denotes the loss probability in the corresponding
$GI/M^{C}/1/N_i$ queueing system, the probability distribution of
interarrival time of which is $B_i(x)$. All probabilities $p_k$,
$1\leq k\leq\ell$ are very small as is $N_k$ large. They decrease
geometrically fast (see Theorem \ref{prop5} below), and the finite
sum of these probabilities seems to remain a good upper bound for
the buffers loss probability.

\begin{thm}
\label{prop5} The buffer contents loss probability is not greater
than $p_1+p_2+\ldots+p_k$, $k=1,2,\ldots,\ell$, where
\begin{equation}
\label{8.1}
\begin{aligned}
p_k=&\frac{(1-\rho_k)[1+C\mu\widehat{B}_k^\prime(\mu-\mu\varsigma_k^C)]
\varsigma_k^{N_k}}
{(1-\rho_k)(1+\varsigma_k+\ldots+\varsigma_k^{C-1})-\rho_k[1+C\mu\widehat{B}_k^\prime(\mu-\mu\varsigma_k^C)]
\varsigma_k^{N_k}}\\
&+o\left(\varsigma_k^{2N_k}\right),
\end{aligned}
\end{equation}
\begin{equation*}
\rho_k=\frac{\lambda_k}{C\mu},
\end{equation*}
and $\varsigma_k$ is the (least) root of the functional equation
$$
z=\widehat{B}_k(\mu-\mu z^C)
$$
in the interval (0,1).
\end{thm}
\begin{proof}
We consider the $GI/M^{C}/1/N_k$ queueing system. Following
Miyazawa \cite{Miyazawa 1990}, the loss probability for the
$GI/M^{Y}/1/N_k$ queueing system is determined by the formula
\begin{equation*}
p_k=\frac{1}{\sum_{j=0}^{N_k}\pi_{k,j}},
\end{equation*}
where the generating function of $\pi_{k,j}$, $j=1,2,\ldots$ is
\begin{equation}\label{8.2}
\Pi_k(z)=\sum_{j=0}^\infty\pi_{k,j}z^j=\frac{(1-Y(z))\widehat{B}_k(\mu-\mu
Y(z))}{\widehat{B}_k(\mu-\mu Y(z))-z},
\end{equation}
and $Y(z)$ is the generating function of complete service batch.
In the case of the $GI/M^{C}/1/N_k$ queueing system $Y(z)=z^C$,
and \eqref{8.2} can be then rewritten as
\begin{equation}\label{8.3}
\Pi_k(z)=\frac{(1-z^C)\widehat{B}_k(\mu-\mu
z^C)}{\widehat{B}_k(\mu-\mu z^C)-z}.
\end{equation}
In the particular case of $C=1$ the asymptotic behaviour of the
loss probability has been studied in \cite{Abramov 2002} and
\cite{Choi Kim Wee 2000}. In the case of $\rho_k<1$ it was based
on an application of the Tak\'acs theorem \cite{Takacs 1967}, p.
22-23.

In the case of $C>1$ the scheme of the proof is similar. Expanding
(1-$z^C$) in the numerator of \eqref{8.3} as $
1-z^C=(1-z)(1+z+\ldots+z^{C-1}) $ we have
\begin{equation}
\label{8.4}\Pi_k(z)=\frac{(1-z)(1+z+z^2+\ldots+z^{C-1})\widehat{B}_k(\mu-\mu
z^C)}{\widehat{B}_k(\mu-\mu z^C)-z}.
\end{equation}
Therefore, the other generating function
$\widetilde{\Pi}_k(z)=\frac{1}{1-z}\Pi_k(z)$ is
\begin{equation}
\label{8.5}
\begin{aligned}
\widetilde{\Pi}_k(z)&=\sum_{i=0}^\infty\widetilde{\pi}_{k,i}z^i
=\sum_{i=0}^\infty\left(\sum_{j=0}^i{\pi}_{k,j}\right)z^i\\&=\frac{(1+z+z^2+\ldots+z^{C-1})\widehat{B}_k(\mu-\mu
z^C)}{\widehat{B}_k(\mu-\mu z^C)-z},
\end{aligned}
\end{equation}
and the loss probability is
\begin{equation}
\label{8.6}p_k=\frac{1}{\widetilde{\pi}_{k,N_k}}.
\end{equation}
Our goal is therefore to find the asymptotic behaviour of
$\widetilde{\pi}_{k,N_k}$ as $N_k\to\infty$.

The equation $z=\widehat{B}_k(\mu-\mu z^C)$ has exactly one
solution $\varsigma_k$ in the interval (0,1). Furthermore,
$\widehat{B}_k(\mu-\mu z^C)$ is the probability generating
function of some integer random variable, i.e.
$$
\widehat{B}_k(\mu-\mu z^C)=R(z)=\sum_{j=0}^\infty r_jz^j,
$$
and
$$
\widetilde{\Pi}_k(z)=F(z)\sum_{i=0}^{C-1}z^i,
$$
where
$$
F(z)=\sum_{n=0}^\infty f_nz^n=\frac{R(z)}{R(z)-z}.
$$
Therefore (see Tak\'acs \cite{Takacs 1967}), the sequence
$\{f_n\}$ satisfies the recurrence relation
$$
f_n=\sum_{j=0}^n f_{n-j+1}r_j, \ \ f_0>0.
$$

Since $\varsigma_k<1$, we correspondingly obtain
$\gamma=\sum_{n=1}^\infty nr_n>1$. According to formula (35) of
\cite{Takacs 1967}, p. 23
\begin{equation}\label{8.14}
\lim_{n\to\infty}\left[f_n-\frac{f_0\delta^{-n}}{[1-F^\prime(\delta)]}\right]=\frac{f_0}{1-\gamma},
\end{equation}
where $\delta$ is the least root of equation $z=R(z)$ in the
interval (0,1).

In our case $\rho_k=\frac{1}{\gamma}$, and $\varsigma_k=\delta$,
and we have:
\begin{equation*}
\begin{aligned}
\label{8.15}\widetilde{\pi}_{k,N_k}=&\frac{1+\varsigma_k+\ldots+\varsigma_k^{C-1}}
{\varsigma_k^{N_k}}\cdot\frac{1}{1+C\mu\widehat{B}_k^\prime(\mu-\mu\varsigma_k^C)}\\
&+\frac{(1+\varsigma_k+\ldots+\varsigma_k^{C-1})\rho_k}{\rho_k-1}+o(1),
\end{aligned}
\end{equation*}
and the statement of the proposition follows from \eqref{8.6}
after some algebraic transformations.
\end{proof}

\subsection{Loss probabilities under heavy load
conditions}\label{sec8.2} The loss probabilities under heavy load
conditions for $GI/M/1/n$ queues have been recently studied in
\cite{Abramov 2002} and \cite{Whitt 2004}. For the further
development of these results see also \cite{Abramov 2006} and
\cite{Whitt 2005}. In this specific case, the behaviour of the
system under heavy load condition differs from the classic cases
considered in these papers.

We consider the case of heavy load conditions and assume that the
load parameter $\rho_\ell=\frac{\lambda_\ell}{C\mu}$ is close to
1. More specifically, we assume that
$\rho_\ell=\rho_\ell(\delta)<$1 ($\delta$ is a small parameter),
and $\rho_\ell(\delta)$ approaches 1 from the left as $\delta$
vanishes. Denote $\rho_{\ell,j}=\mu^j\int_0^\infty
x^j\mathrm{d}B_\ell(x)$ \
($\rho_{\ell,1}=\frac{1}{C\rho_{\ell}})$. We have the following
result.

\begin{thm}
\label{prop6} Assume that $\rho_\ell(\delta)<1$,
$\rho_\ell(\delta)$ approaches 1 from the left, and $\delta
N_\ell(\delta)\to \Delta>0$ as $\delta$ vanishes. Assume also that
$\widetilde\rho_{\ell,2}=\lim_{\delta\to 0}\rho_{\ell,2}(\delta)$,
and $\rho_{\ell,3}(\delta)$ remains bounded as $\delta$ vanishes.
Then
\begin{equation}\label{8.18}
p_\ell=\frac{\delta\exp\left(-\dfrac{\Delta}{\binom{C}{2}\widetilde\rho_{\ell,2}}\right)}
{C-\exp\left(-\dfrac{\Delta}{\binom{C}{2}\widetilde\rho_{\ell,2}}
\right)} \cdot[1+o(1)].
\end{equation}
\end{thm}

\begin{proof}
Let us first derive the expansion for the least root of equation
$z=\widehat B_\ell(\mu-\mu z^C)$ under the assumption of the
theorem. Clearly, the root of this equation approaches 1 as
$\delta$ vanishes. Therefore, using the Taylor expansion of
$\widehat B_\ell(\mu-\mu z^C)$ as $\delta$ vanishes, we obtain the
following equation for $\varsigma_\ell$
\begin{equation*}
\label{8.19}
z=1-(1-\delta)(1-z)+\binom{C}{2}\widetilde\rho_{\ell,2}(1-z)^2+o(1-z)^3.
\end{equation*}
Ignoring the last term $o(1-z)^3$ we have the quadratic equation,
the solutions of which are $z=1$ and
$z=1-\delta/\left[\binom{C}{2}\widetilde\rho_{\ell,2}\right]$.
Therefore we obtain
\begin{equation}
\label{8.20}
\varsigma_\ell=1-\frac{\delta}{\binom{C}{2}\widetilde\rho_{\ell,2}}+o(\delta).
\end{equation}
Notice, that representation similar to \eqref{8.20} for the root
of equation $z=\widehat B_\ell(\mu-\mu z)$ (particular case where
$C$=1) has been obtained in Subhankulov \cite{Subhankulov 1976},
p. 326.

Next, the asymptotic representation for $p_k$, $k=1,2,\ldots,\ell$
is given by \eqref{8.1}. For $k=\ell$ the main term of asymptotic
expansion of
$$
[1+C\mu\widehat{B}_\ell^\prime(\mu-\mu\varsigma_\ell^C)]
\varsigma_\ell^{N_\ell}
$$
is given by
\begin{equation*}
\delta\exp\left({-\frac{\Delta}{\binom{C}{2}\widetilde\rho_{\ell,2}}}\right)
\end{equation*}
(see \cite{Abramov 2002}, \cite{Abramov 2006} for details of the
proof), and according to \eqref{8.20} the main term of the
asymptotic expansion of
$$
1+\varsigma_\ell+\ldots+\varsigma_\ell^{C-1}
$$
 is given
by
\begin{equation*}
C-\frac{ \delta}{\widetilde\rho_{\ell,2}}.
\end{equation*}
Next notice, that the expansion for the term
\begin{equation*}
1-\rho_\ell-\rho_\ell[1+C\mu\widehat{B}_k^\prime(\mu-\mu\varsigma_\ell^C)]
\varsigma_\ell^{N_\ell}
\end{equation*}
is
\begin{equation*}
\delta\left[1-\exp\left({-\frac{\Delta}
{\binom{C}{2}\widetilde\rho_{\ell,2}}}\right)\right][1+o(1)].
\end{equation*}
Therefore, asymptotic relation \eqref{8.18} follows.
\end{proof}

\section{Approximation of the solution in particular
cases}\label{sec9} In this section we discuss the approximation of
the solution for the problem stated in Section \ref{sec1.4} that
is to minimize functional \eqref{1.4+} containing the terms
$\alpha^{(k)}J^{(k)}$ associated with buffer contents
$\mathcal{Q}^{(k)}(t)$, $k=1,2,\ldots,\ell$.

However, in all particular cases above the explicit solutions were
obtained for the cumulative buffer contents $\mathcal{Q}_k(t)$,
$k=1,2,\ldots,\ell$, and the solution of the problem in the
initial terms seems to be hard. Therefore we formulate and solve
the problem in new terms. This solution of the new problem is then
used to approximate the desired solution of the initial problem.

Let us first introduce new functionals $J_k$ instead of the
$J^{(k)}$ ($k=1,2$,\ldots, $\ell$), which were introduced in
Section \ref{sec1.4}.

Namely, let $N_1$, $N_2$,\ldots,$N_\ell$ denote large integer
numbers, $N_1<N_2<\ldots<N_\ell$. We set
\begin{eqnarray*}\label{6.0}
J_{k}=\mathbb{P}^{\_}\lim_{t\to\infty}\frac{1}{{A}_\ell(t)}
\sum_{j=1}^{{A}_{k}(t)}\mathbf{I}\left\{\mathcal{Q}_{k}(t_{k,j})>
N_{k}\right\},\\
k=1,2,\ldots,\ell.\nonumber
\end{eqnarray*}

Note first (see \eqref{4.9}) that for $J_k$, $k=1,2,\ldots,\ell$
we have the following representation:
\begin{equation}
\label{6.0add}
\begin{aligned}
J_k&=\frac{1}{\lambda_\ell}\lim_{t\to\infty}\frac{1}{t}\mathbb{E}\sum_{l=1}^{C}
\int_0^t\mathbf{I}
\left\{\mathcal{Q}_k(u-)\geq N_{k}+l\right\}\mathrm{d}{D}(u),\\
&\ \ \ \ \ \ k=1,2,\ldots,\ell,
\end{aligned}
\end{equation}
where
$\lambda_k=\lambda^{(1)}+\lambda^{(2)}+\ldots+\lambda^{(k)}$,
$k=1,2,\ldots,\ell$. The proof of representation \eqref{6.0add} is
similar to the proof of \eqref{4.9} with minor difference in the
notation.

Replacing functional \eqref{1.4+}  by
\begin{equation}\label{9.1}
\overline{J}=\alpha_1J_1+\alpha_2J_2+\ldots+\alpha_\ell J_\ell,
\end{equation}
we have then the following problems similar  to the problems
formulated above in Section \ref{sec1.4}.

1. Assuming that $N_1$, $N_2$,\ldots,$N_\ell$ are known, minimize
$C$ subject to $\overline{J}\leq\varepsilon$.

2. Assume that $C$ is known, but $N_1$, $N_2$,\ldots,$N_\ell$ are
unknown. Assume additionally that with given $\beta_1$,
$\beta_2$,\ldots,$\beta_{\ell-1}$ the values $N_1$, $N_2$,\ldots,
$N_\ell$ must satisfy the condition: $N_1=\lfloor\beta_1
N_2\rfloor$=$\lfloor\beta_2 N_3\rfloor$=\ldots
=$\lfloor\beta_{\ell-1} N_\ell\rfloor$. The problem is to minimize
$N_1$ subject to $\overline{J}\leq\varepsilon$.

The values $\alpha_1$, $\alpha_2$, \ldots, $\alpha_\ell$; $N_1$,
$N_2$,\ldots, $N_\ell$ are unknown, and by approximation of the
solution of the problem we hope to find a correspondence between
the vectors ($\alpha_1$, $\alpha_2$, \ldots, $\alpha_\ell$) and
($\alpha^{(1)}$, $\alpha^{(2)}$, \ldots, $\alpha^{(\ell)}$) and
between the vectors ($N_1$, $N_2$,\ldots, $N_\ell$) and
($N^{(1)}$, $N^{(2)}$,\ldots, $N^{(\ell)}$) such that the solution
of the initial problems formulated in Section \ref{sec1.4} and the
problems formulated in this section would be approximately the
same.

Consider first the queueing systems with infinite number of
waiting places, say $GI/M^C/1$ queues.

Notice, that $\alpha_1=\alpha^{(1)}$ and $N_1=N^{(1)}$. According
to Theorem \ref{prop3} the expected queue-length of the $k$th
cumulative buffer content immediately before arrival of a unit is
\begin{equation}\label{9.2}
\begin{aligned}
\lim_{t\to\infty}\frac{1}{t}\sum_{m=1}^\infty
m\int_0^t\mathbf{I}\{\mathcal{Q}_k(u)=m\}
\mathrm{d}u=\frac{\varsigma_k}{1-\varsigma_k}.
\end{aligned}
\end{equation}
From \eqref{9.2} we have the following. Put
\begin{equation}\label{9.3}
\begin{aligned}
&p_{2,1}=\frac{\dfrac{\varsigma_1}{1-\varsigma_1}}
{\dfrac{\varsigma_2}{1-\varsigma_2}}=\frac{\varsigma_1(1-\varsigma_2)}{\varsigma_2(1-\varsigma_1)},\\
&p_{2,2}=1-p_{2,1},
\end{aligned}
\end{equation}
and then
\begin{equation}
\label{9.4} \alpha_2=\alpha_1p_{2,1}+\alpha^{(2)}p_{2,2}.
\end{equation}
Similarly to \eqref{9.3} and \eqref{9.4} for $k=1,2,\ldots,\ell-1$
we set
\begin{equation}\label{9.5}
\begin{aligned}
&p_{k+1,1}=\frac{\varsigma_k(1-\varsigma_{k+1})}{\varsigma_{k+1}(1-\varsigma_{k})},\\
&p_{k+1,2}=1-p_{k+1,1},
\end{aligned}
\end{equation}
and
\begin{equation}
\label{9.6} \alpha_{k+1}=\alpha_kp_{k+1,1}
+\alpha^{(k+1)}p_{k+1,2}.
\end{equation}
Let us now express the correspondence between the vectors
$(N^{(1)}$, $N^{(2)}$,\ldots, $N^{(\ell)}$) and $(N_1$,
$N_2$,\ldots, $N_\ell$). Let $\beta^{(2)}$, $\beta^{(3)}$,\ldots,
$\beta^{(\ell)}$ be such the real numbers that
$N^{(1)}=\lfloor\beta^{(2)} N^{(2)}\rfloor$=$\lfloor\beta^{(3)}
N^{(3)}\rfloor$ = \ldots = $\lfloor\beta^{(\ell)}
N^{(\ell)}\rfloor$.

Then for the purpose of approximation the values $\beta_2$,
$\beta_3$,\ldots,$\beta_{\ell}$ are taken as
\begin{equation}
\label{9.7} \beta_{k}=\frac{\beta^{(k)}}{1+\beta^{(k)}}, \
k=2,3\ldots,\ell,
\end{equation}
and $N_1=\lfloor\beta_2 N_2\rfloor$=$\lfloor\beta_3 N_3\rfloor$ =
\ldots = $\lfloor\beta_\ell N_\ell\rfloor$.

For the queueing model with large finite buffers, say
$GI/M^C/1/N_k$, $k=1,2,\ldots,\ell$, the approximation is similar.
Specifically, since the buffers are large, approximation for
($\alpha_1$, $\alpha_2$,\ldots, $\alpha_\ell$) can be given by
\eqref{9.2}-\eqref{9.5}. The values $\beta_2$,
$\beta_3$,\ldots,$\beta_\ell$ are assumed to be taken by the same
relation \eqref{9.7}.

\section{Minimization algorithms for the functional $\overline{J}$}
\label{sec10} In this section we discuss the problem of
minimization of the functional $\overline{J}$ defined by
\eqref{9.1}. $J_1$, $J_2$,\ldots, $J_\ell$ depends on parameters
$C$, $N_1$, $N_2$,\ldots, $N_\ell$.

\subsection{$C$ is known while $N_1$ is unknown}
Assume first that $C$ is known, $N_1$ is unknown, and the problem
is to find the value $N_1$ minimizing the functional
$\overline{J}$ in the buffer models, where explicit representation
for the state probabilities as well as for $J_k$,
$k=1,2,\ldots,\ell$ are known. These models are considered in
Sections  \ref{sec7} and \ref{sec8}.

\smallskip
 To be specific we refer to
the models of infinite buffers of $GI/M^C/1$ queues. The algorithm
has the following steps.

\smallskip
$\bullet$ \textit{Step 1}. Calculate $\varsigma_k$,
$k=1,2,\ldots,\ell$. Recall that $\varsigma_k$ is the root of the
functional equation \eqref{7.3} in the interval (0,1). For each
$k$ it can be calculated by the fixed point method or by one of
other well-known methods, say direct search method or gold section
method (e.g. see \cite{Himmelblau 1972}).

\smallskip
$\bullet$ \textit{Step 2}. We have $\ell$ geometric distributions
obtained in Step 1, and therefore one can compute the
corresponding values $N_1$, $N_2$,\ldots,$N_\ell$ at which each of
the tails of the geometric distributions multiplied to the
corresponding coefficient $\alpha_k$, $k=1,2,\ldots,\ell$ will be
less than $\varepsilon$ (i.e. $\alpha_kJ_k<\varepsilon$).

\smallskip
$\bullet$ \textit{Step 3}. By using the known coefficients
$\beta_2$, $\beta_3$, \ldots, $\beta_{\ell}$ one can find the
value $N_1^{lower}$. $N_1^{lower}$ is the maximum amongst all
minimal values of $N_1$ under which $\alpha_kJ_k<\varepsilon$ for
all $k=1,2,\ldots,\ell$. Specifically, we have the system:
\begin{equation*}
\begin{aligned}
&\min N_1: \ \alpha_1J_1(N_1)<\varepsilon,\\
&\min N_1: \ \alpha_2J_2(N_1)<\varepsilon,\\
&......................................\\
&\min N_1: \ \alpha_\ell J_\ell(N_1)<\varepsilon,
\end{aligned}
\end{equation*}
and $N_1^{lower}$ is the maximum amongst $\ell$ obtained values of
$N_1$.

\smallskip
$\bullet$ \textit{Step 4}. By using the same known coefficients
$\beta_2$, $\beta_3$, \ldots, $\beta_{\ell}$ one can find the
value $N_1^{upper}$. $N_1^{upper}$ is the minimum amongst all
maximal values of $N_1$ under which
$\alpha_kJ_k<\frac{\varepsilon}{\ell}$ for all
$k=1,2,\ldots,\ell$. Specifically, we have the system:
\begin{equation*}
\begin{aligned}
&\min N_1: \ \alpha_1J_1(N_1)<\frac{\varepsilon}{\ell},\\
&\min N_1: \ \alpha_2J_2(N_1)<\frac{\varepsilon}{\ell},\\
&......................................\\
&\min N_1: \ \alpha_\ell J_\ell(N_1)<\frac{\varepsilon}{\ell},
\end{aligned}
\end{equation*}
and $N_1^{upper}$ is the maximum amongst $\ell$ obtained values of
$N_1$.

\smallskip
$\bullet$ \textit{Step 5}. We solve  the following integer
programming problem:
\begin{equation*}
\mbox{minimize} \ N_1: N_1^{lower}\leq N_1\leq N_1^{upper}
\end{equation*}
\begin{equation*}
\mbox{subject to} \ \alpha_1J_1+\alpha_2J_2+\ldots+\alpha_\ell
J_{\ell}\leq\varepsilon.
\end{equation*}

\subsection{$N_1$, $N_2$,\ldots,$N_\ell$ are known while $C$ is
unknown} In the case where $N_1$, $N_2$,\ldots, $N_\ell$ all are
known but $C$ is unknown the algorithm of the problem solution is
the following.
\smallskip

$\bullet$ \textit{Step 1}. From the stability condition find the
lower (integer) bound for $C$:
\begin{equation*}
C^{lower}=\min\left\{C: \frac{\lambda_\ell}{C\mu}<1\right\}.
\end{equation*}

\smallskip
$\bullet$ \textit{Step 2}. Find $\varsigma_k$,
$k=1,2,\ldots,\ell$.

\smallskip
$\bullet$ \textit{Step 3}. Compute the functional $\overline{J}$.

\smallskip
$\bullet$ If $\overline{J}>\varepsilon$, then find a new value $C$
and repeat steps 1-3. These procedure should be repeated more and
more while $\overline{J}>\varepsilon$. Since the upper bound of
$C$ is unknown, the value $C$ should be found according to the
special search procedure offered by Rubalskii \cite{Rubalskii
1982}.

Rubalskii \cite{Rubalskii 1982} proposed the minimization
algorithm for a unimodal function on an unbounded set. The optimal
algorithm is an extension of the standard Fibonacci procedure.

\section{Concluding remarks}\label{sec12}
In this paper we studied queueing systems with priority classes
and infinite and finite buffers. We derived general type equations
for buffer content processes assuming that service mechanism is
autonomous. The results of general theory were then applied to
special queueing models with exponentially distributed service
times. These queueing systems are a particular case of systems
with an autonomous service mechanism. For the model having large
buffers we derived an asymptotic result for the loss probability.
We developed an algorithm for a solution of the problem
numerically.

\subsection*{Acknowledgements}
The author expresses his gratitude to Prof. Robert Liptser for
providing the author relevant information about this problem. The
author indebts to Prof. Boris Miller for his advice and some
comments. Careful reading and comments of Prof. Aidan Sudbury are
appreciated. The author also thanks the anonymous referee for
careful reading and comments, which substantially improved the
presentation. The research was supported by the Australian
Research Council, grant \# DP0771338.

\section*{APPENDIX A: Proof of Lemma \ref{lem1}}
We start from equation \eqref{1.1}. In order to write this
equation in the customary form of a stochastic equation, we use
the process $D(t)=\sum_{n=1}^\infty \mathbf{I}\{x_n\leq t\}$. The
jumps of the process $D(t)$ are equal to 1, and according to the
definition, we have $CD(t)=\mathcal{D}(t)$ for all $t\geq 0$. Then
due to the assumption that arrival and departure jumps are
disjoint, \eqref{1.1} can be rewritten
\begin{equation*}
\mathcal{Q}^{(1)}(t)=\mathcal{A}^{(1)}(t)-\sum_{j=1}^C\int_0^t
\mathbf{I}{\left\{\mathcal{Q}^{(1)}(u-)\geq j\right\}}
\mathrm{d}D(u).\leqno(A.1)
\end{equation*}
The equivalence of representations of \eqref{1.1} and (A.1) can be
easily checked by considering a small time interval ($t-\delta,
t$] containing exactly one event as either arrival or departure of
a unit. Then the term
\begin{equation*}
\mathbf{I}{\left\{\mathcal{Q}^{(1)}(u-)\geq j\right\}}\leqno(A.2)
\end{equation*}
of the integrand shows that if $u$ is the point of jump of the
process $\mathcal{D}(t)$, and $\mathcal{Q}^{(1)}(u-)=n\leq C$,
then $\mathcal{Q}^{(1)}(u)=0$.

Similarly to (A.1), equation \eqref{1.2} can be rewritten as
follows:
$$
\mathcal{Q}^{(2)}(t)=\mathcal{A}^{(2)}(t)-\sum_{j=1}^C\int_0^t
\mathbf{I}\left\{\mathcal{Q}^{(1)}(u)=0\right\}\leqno(A.3)
$$
$$
\times\mathbf{I}\left\{\mathcal{Q}^{(2)}(u-)\geq j
-\mathcal{Q}^{(1)}(u-)\right\} \mathrm{d}D(u).
$$
The explanation of the equivalence of (A.1) and (A.3) is similar
to the above case, but slightly more complicated in details.
Specifically, the presence of the term
$\mathbf{I}\left\{\mathcal{Q}^{(1)}(u)=0\right\}$ in the integrand
is obvious, and the validation of the term
\begin{equation*}
\mathbf{I}\left\{\mathcal{Q}^{(2)}(u-) \geq j
-\mathcal{Q}^{(1)}(u-)\right\}
\end{equation*}
is explained similarly to that of (A.2).

Let us now find the representation for $\mathcal{Q}_2(t)$ =
$\mathcal{Q}^{(1)}(t)$ + $\mathcal{Q}^{(2)}(t)$. Keeping in mind
that $\mathcal{Q}_1(t)=\mathcal{Q}^{(1)}(t)$ and
$\mathcal{A}_1(t)=\mathcal{A}^{(1)}(t)$ from (A.1) and (A.3) we
obtain:
\begin{eqnarray*}
\mathcal{Q}_2(t) &=&\mathcal{A}_2(t)
-\sum_{j=1}^C\int_0^t\mathbf{I}{\{\mathcal{Q}_1(u-)\geq j\}}
 \mathrm{d}D(u)\nonumber\\
&&-\sum_{j=1}^C\int_0^t\mathbf{I}{\{\mathcal{Q}_1(u-)<j\}}
\mathbf{I}\{\mathcal{Q}_2(u-)\geq j\}\mathrm{d}D(u)\nonumber\\
&=&\mathcal{A}_2(t)-\sum_{j=1}^C\int_0^t\mathbf{I}\{\mathcal{Q}_2(u-)\geq
j\} \mathrm{d}D(u).
\end{eqnarray*}
The term
$$
\sum_{j=1}^C\int_0^t\mathbf{I}{\{\mathcal{Q}_1(u-)\geq j\}}
 \mathrm{d}D(u)
$$
characterizes departure lengths from the highest priority buffer,
while the term
$$
\sum_{j=1}^C\int_0^t\mathbf{I}{\{\mathcal{Q}_1(u-)<j\}}
\mathbf{I}\{\mathcal{Q}_2(u-)\geq j\}\mathrm{d}D(u)
$$
characterizes that from the second buffer of the lower priority.

Thus for $k=1,2$ we have already shown
\begin{equation*}
\begin{aligned}
\mathcal{Q}_k(t)=\mathcal{A}_k(t)-\sum_{j=1}^C\int_0^t\mathbf{I}\{\mathcal{Q}_k(u-)
\geq j\} \mathrm{d}D(u).
\end{aligned}
\end{equation*}
Let us prove \eqref{2.0} by using induction.

For this purpose let us write first a representation for
$\mathcal{Q}^{(k+1)}(t)$, $k=1,2,\ldots,\ell-1$. Similarly to
(A.3) we have:
\begin{equation*}
\mathcal{Q}^{(k+1)}(t)=\mathcal{A}^{(k+1)}(t)-\sum_{j=1}^C\int_0^t
\mathbf{I}\left\{\mathcal{Q}_k(u)=0\right\}\leqno(A.4)
\end{equation*}
\begin{equation*}
\times\mathbf{I}\left\{\mathcal{Q}^{(k+1)}(u-)\geq j
-\mathcal{Q}_k(u-)\right\} \mathrm{d}D(u).
\end{equation*}
Equation (A.4) is a straightforward extension of (A.3). Therefore,
assuming that \eqref{2.0} is valid for some $k$ and adding
$\mathcal{Q}_k(t)$ and $\mathcal{Q}^{(k+1)}(t)$, and similarly to
the above for the $k+1$st cumulative  buffer content we obtain:
\begin{eqnarray*}
\mathcal{Q}_{k+1}(t)
&=&\mathcal{A}_{k+1}(t)-\sum_{j=1}^C\int_0^t\mathbf{I}\{\mathcal{Q}_{k+1}(u-)
\geq j\} \mathrm{d}D(u).
\end{eqnarray*}
Representation \eqref{2.0} is proved.

\section*{APPENDIX B: Deriving \eqref{4.9}}
Equation \eqref{4.4} is a basic equation for our analysis. With
$\infty\cdot0=0$ for the left-hand side of \eqref{4.4} we have:
\begin{equation*}
\begin{aligned}
&\lim_{t\to\infty}\frac{1}{t}\mathbb{E}\sum_{i=1}^{A^{(k)}(t)}\mathbf{I}\left\{\mathcal{Q}^{(k)}\Big(t_{i}^{(k)}\Big)\geq
m, \ \mathcal{Q}^{(k)}\Big(t_{i}^{(k)}-\Big)<m\right\}\\
&=\lambda^{(k)}\lim_{t\to\infty}
\mathbb{E}\frac{1}{A^{(k)}(t)}\sum_{i=1}^{A^{(k)}(t)}\mathbf{I}\left\{\mathcal{Q}^{(k)}\Big(t_{i}^{(k)}\Big)\geq
m, \ \mathcal{Q}^{(k)}\Big(t_{i}^{(k)}-\Big)<m\right\}\\
&=\lambda^{(k)}\lim_{t\to\infty}\sum_{l=0}^\infty\frac{1}{l}
\sum_{i=1}^{l}\mathbb{P}\left\{\mathcal{Q}^{(k)}\Big(t_{i}^{(k)}\Big)\geq
m, \ \mathcal{Q}^{(k)}\Big(t_{i}^{(k)}-\Big)<m\right\}\\
&\ \ \ \times\mathbb{P}\left\{A^{(k)}(t)=l\right\}\\
&=\lambda^{(k)}\lim_{t\to\infty}\mathbb{E}\sum_{l=0}^\infty\frac{1}{l}
\sum_{i=1}^{l}\mathbf{I}\left\{\mathcal{Q}^{(k)}\Big(t_{i}^{(k)}\Big)\geq
m, \ \mathcal{Q}^{(k)}\Big(t_{i}^{(k)}-\Big)<m\right\}\\
&\ \ \ \times\mathbb{P}\left\{A^{(k)}(t)=l\right\}\\
&=\lambda^{(k)}\mathbb{P}^{\_}\lim_{t\to\infty}
\frac{1}{A^{(k)}(t)}\sum_{i=1}^{A^{(k)}(t)}\mathbf{I}\left\{\mathcal{Q}^{(k)}\Big(t_{i}^{(k)}\Big)\geq
m, \ \mathcal{Q}^{(k)}\Big(t_{i}^{(k)}-\Big)<m\right\},\\
&\ \ \ \ \ \ k=1,2,\ldots,\ell.
\end{aligned}
\end{equation*}
The right-hand side of \eqref{4.4} can be rewritten as follows:
\begin{equation*}
\begin{aligned}
&\lim_{t\to\infty}\frac{1}{t}\mathbb{E}\sum_{l=1}^C\sum_{j=1}^{D(t)}\mathbf{I}
\left\{\mathcal{Q}^{(k)}(x_{j}-)=m-1+l\right\}\\
&=\lim_{t\to\infty}\frac{1}{t}\mathbb{E}\sum_{l=1}^C\int_0^t\mathbf{I}
\left\{\mathcal{Q}^{(k)}(u-)=m-1+l\right\}\mathrm{d}D(u),\\
 &\ \ \ \ \ \ k=1,2,\ldots,\ell.
\end{aligned}
\end{equation*}
From these last two equations we obtain:
\begin{equation*}
\begin{aligned}
&\mathbb{P}^{\_}\lim_{t\to\infty}
\frac{1}{A^{(k)}(t)}\sum_{i=1}^{A^{(k)}(t)}\mathbf{I}\left\{\mathcal{Q}^{(k)}\Big(
t_{i}^{(k)}\Big)\geq
m, \ \mathcal{Q}^{(k)}\Big(t_{i}^{(k)}-\Big)<m\right\}\\
&=\frac{1}{\lambda^{(k)}}\lim_{t\to\infty}\frac{1}{t}\mathbb{E}
\sum_{l=1}^C\int_0^t\mathbf{I}
\left\{\mathcal{Q}^{(k)}(u-)=m-1+l\right\}\mathrm{d}D(u),\\
&\ \ \ \ \ \ k=1,2,\ldots,\ell,
\end{aligned}
\end{equation*}
and taking into account \eqref{1.4} we finally obtain \eqref{4.9}.


\begin{thebibliography}{10}
\bibitem{Abramov 2000}\textsc{Abramov, V.M.} (2000). A large
closed queueing network with autonomous service and bottleneck.
\emph{Queueing Systems}, 35: 23-54.

\smallskip

\bibitem{Abramov 2002}\textsc{Abramov, V.M.} (2002). Asymptotic
analysis of the $GI/M/1/n$ loss system as $n$ increases to
infinity. \emph{Annals of Operations Research}, 112: 35-41.

\smallskip

\bibitem{Abramov 2004}\textsc{Abramov, V.M.} (2004). Large closed
queueing networks containing two types of node and multiple
customer classes: One bottleneck station. \emph{Queueing Systems},
48: 45-73.

\smallskip

\bibitem{Abramov 2004 SIAM}\textsc{Abramov, V.M.} (2004).
Asymptotic behaviour of the number of lost messages. \emph{SIAM
Journal on Applied Mathematics}, 64: 746-761.

\smallskip

\bibitem{Abramov 2005}\textsc{Abramov, V.M.} (2005). The stability
of join-the-shortest-queue models with general input and output
processes. arXiv : math/PR 0505040.

\smallskip

\bibitem{Abramov 2006} \textsc{Abramov, V.M.} (2007). Asymptotic
analysis of loss probabilities in $GI/M/m/n$ queueing systems as
$n$ increases to infinity. \emph{Quality Technology and
Quantitative Management}, 4: 379-393.

\smallskip
\bibitem{Abramov 2008} \textsc{Abramov, V.M.} (2008). Large closed
queueing networks in semi-Markov environment and their
applications. \emph{Acta Applicandae Mathematicae}, 100: 201-226.

\smallskip
\bibitem{Anulova and Liptser 1990} \textsc{Anulova, S.V. and Liptser,
R. Sh.} (1990).  Diffusion approximation for the processes with
normal reflection. \emph{Theory of Probability and Its
Application}, 35: 413-423.

\smallskip

\bibitem{Berger and Whitt 1998a}\textsc{Berger, A.W. and Whitt,
W.} (1998). Effective bandwidths with priorities. \emph{IEEE/ACM
Transaction on Networking}, 6 (4): 447-460.

\smallskip

\bibitem{Berger and Whitt 1998b}\textsc{Berger, A.W. and Whitt,
W.} (1998). Extending the effective bandwidth concept to networks
with priority classes. \emph{IEEE Communication Magazine}, 36 (8):
78-83.

\smallskip

\bibitem{Bestimas et al 1998}\textsc{Bertsimas, D., Paschalidis,
I.C. and Tsitsiklis, J.N.} (1998). Asymptotic buffer overflow
probabilities in multiclass multiplexers: An optimal control
approach. \emph{IEEE Transactions on Automatic Control} 43:
315-335.


\smallskip


\bibitem{Borovkov 1976}\textsc{Borovkov, A.A.} (1976). \emph{Stochastic
Processes in Queueing Theory}. Springer, Berlin.

\smallskip

\bibitem{Borovkov 1984}\textsc{Borovkov, A.A.} (1984). \emph{Asymptotic
Methods in Queueing Theory}. John Wiley, New York.

\smallskip

\bibitem{Botvich and Duffield 1995}\textsc{Botvich, D.D. and
Duffield, N.G.} (1995). Large deviations, the shape of loss curve,
and economics of scale in large multiplexers. \emph{Queueing
Systems}: 20, 293-320.


\smallskip


\bibitem{Chao Pinedo and Shaw 1996}\textsc{Chao, X., Pinedo, M.
and Shaw, D.} (1996). Network of queues with batch services and
customer coalescence. \emph{Journal of Applied Probability} 33:
858-869.

\smallskip
\bibitem{Choi and Kim (2000)} \textsc{Choi, B.D. and Kim, B.} (2000). Sharp
results on convergence rates for the distribution of the
$GI/M/1/K$ queues as $K$ tends to infinity. \emph{Journal of
Applied Probability} 37: 1010-1019.


\smallskip
\bibitem{Choi Kim Wee 2000} \textsc{Choi, B.D., Kim, B. and Wee,
I.-S.} (2000). Asymptotic behavior of loss probability in
$GI/M/1/K$ queue as $K$ tends to infinity. \emph{Queueing
Systems}, 36: 437-442.

\smallskip
\bibitem{Choudhury Lucantoni and Whitt 1996} \textsc{Choudhury,
G.L., Lucantoni, D.M. and Whitt, W.} (1996). Squeezing the most
out of ATM. \emph{IEEE Transactions in Communications} 44:
203-217.



\smallskip



\bibitem{Courcoubetis Siris Stamoulis 1999}\textsc{Courcoubetis, C.,
Siris, V.A. and Stamoulis, G.} (1999). Application of many sources
asymptotic and effective bandwidth for traffic engineering.
\emph{Telecommunication Systems}, 12: 167-191.

\smallskip

\bibitem{Economou and Fakinos 2003}\textsc{Economou, A. and
Fakinos, D.} (2003). On the stationary distribution of the $GI^X
/M^Y/1$ queueing system. \emph{Stochastic Analysis and
Applications}, 21: 559-565.

\smallskip

\bibitem{Elwalid and Mitra 1995}\textsc{Elwalid, A.I. and Mitra, D.}
(1995). Analysis, approximations and admission control of a
multi-service multiplexing system with priorities. \emph{Proc.
IEEE INFOCOM'95}, 463-472.

\smallskip



\bibitem{Elwalid and Mitra 1999}\textsc{Elwalid, A.I. and Mitra, D.}
(1999). Design of generalized processor sharing schedulers with
statistically multiplex heterogeneous QoS classes. \emph{Proc.
IEEE INFOCOM'99}, 1220-1230.

\smallskip

\bibitem{Evans and Everitt
1999}\textsc{Evans, J.S. and Everitt, D.} (1999). Effective
bandwidth-based admission control for multiservice CDMA cellular
networks. \emph{IEEE Transactions on Vehicular Technology}, 48:
36-46.

\smallskip
\bibitem{Filippov 1988} \textsc{Filippov, A.F.} (1988).
\emph{Differential Equations with Discontinuous Right-Hand Side}.
Kluwer, Dordrecht.

\bibitem{Fricker 1986}\textsc{Fricker, C.} (1986).  Etude d'une file GI/G/1
\'a service autonome (avec vacances du serveur).  \emph{Advances
in Applied Probability}, 18: 283-286.

\bibitem{Fricker 1987}\textsc{Fricker, C.} (1987). Note sur un modele de file GI/G/1
\'a service autonom\'e (avec vacances du serveur).  \emph{Advances
in Applied Probability}, 19: 289-291.

\bibitem{Gelenbe and Iasnogorodski 1979}\textsc{Gelenbe, E. and
Iasnogorodski, R.} (1980). A queue with server of walking type
(autonomous service). \emph{Ann. Inst. H. Poincare}, 16, 63-73.

\bibitem{Gnedenko Kovalenko 1968}\textsc{Gnedenko, B.V. and
Kovalenko, I.N.} (1968). \emph{Introduction to the Theory of
Queues}. Israel Program for Scientific Translations, Jerusalem.

\smallskip

\bibitem{Himmelblau 1972}\textsc{Himmelblau, D.M.} (1972).
\emph{Applied Non-Linear Programming}. McGraw-Hill, New York.

\smallskip

\bibitem{Kella 1993}\textsc{Kella, O.} (1993). Parallel and tandem
fluid networks with dependent L\'evy inputs. \emph{The Annals of
Applied Probability}, 3: 682-695.

\smallskip

\bibitem{Kella and Whitt 1992}\textsc{Kella, O. and Whitt, W.}
(1992). A tandem fluid network with L\'evy input. In \emph{Queues
and Related Models} (I. Basawa and U. Bhat eds) 112-128, Oxford
University Press.

\smallskip

\bibitem{Kelly 1996}\textsc{Kelly, F.P.} (1996). Notes on
effective bandwidth. In \emph{Stochastic Networks: Theory and
Applications} (F.P.Kelly, S.Zachary and I.B.Ziedins eds). Oxford
University Press, Oxford, 1996, 141-168.

\smallskip

\bibitem{Kogan and Liptser 1993}\textsc{Kogan, Ya. and Liptser, R. Sh.}
(1993). Limit non-stationary behaviour of large closed queueing
network with bottlenecks. \emph{Queueing Systems}, 14: 33-55.



\smallskip

\bibitem{Kumaran et al 2000}\textsc{Kumaran, K., Margrave, G.E., Mitra, D. and Stanley, K.R.}
(2000). Novel techniques for the design and control of generalized
processor-sharing schedulers for multiple QoS classes. \emph{Proc.
INFOCOM'00}, 2: 932-941.

\smallskip

\bibitem{Lee et al 2005}\textsc{Lee, J.Y., Kim, S., Kim, D. and Sung,
D.K.}(2005). Bandwidth optimization for internet traffic in
generalized processor-sharing servers. \emph{IEEE Transactions on
Parallel Distributed Systems}, 16: 324-334.

\smallskip



\bibitem{Miyazawa 1990}\textsc{Miyazawa, M.} (1990). Complementary
generating functions for the $M^X/GI/1/k$ and $GI/M^Y/1/k$ queues
and their application to the comparison for loss probabilities.
\emph{Journal of Applied Probability}, 27: 684-692.

\smallskip

\bibitem{Paschalidis 1996}\textsc{Paschalidis, I.C.} (1996). Large
deviations in high-speed communications networks. PhD thesis. MIT
Laboratory for Information and Decision Systems, Cambridge, MA,
USA.

\smallskip

\bibitem{Ramanan 2006}\textsc{Ramanan, K.} (2006). Reflected
diffusions defined via extended Skorokhod map. \emph{Electronic
Journal of Probability}, 11: 934-992.

\smallskip

\bibitem{Rubalskii 1982}\textsc{Rubalskii, G.B.} (1982). The
search of an extremum of  unimodal function of one variable in an
unbounded set. \emph{U.S.S.R. Comput. Maths. Math. Phys.}, 22 (1):
8-15. Transl. from  Russian: \emph{Zhurnal Vychislitelnoi
Matematiki i Matematicheskoi Fiziki}, 22 (1): 10-16, 251.

\smallskip
\bibitem{Skorokhod 1961}\textsc{Skorokhod, A.V.} (1961).
Stochastic equations for difusion processes in a bounded region.
\emph{Theory Probabilities and Its Application} 6: 264-274.

\smallskip
\bibitem{Subhankulov 1976} \textsc{Subhankulov, M.A.} (1976). {\em Tauberian
Theorems with Remainder.} Nauka, Moscow. (In Russian.)



\smallskip

\bibitem{Takacs 1962}\textsc{Tak\'acs, L.} (1962). \emph{Introduction to the
Theory of Queues}. Oxford University Press, New York/London.

\smallskip

\bibitem{Takacs 1967}\textsc{Tak\'acs, L.} (1967). \emph{Combinatorial
Methods in the Theory of Stochastic Processes}. John Wiley, New
York.

\smallskip
\bibitem{Tanaka 1979}\textsc{Tanaka, H.} (1979). Stochastic differential equations
with reflected boundary condition in convex regions.
\emph{Hiroshima Mathematical Journal}, 9: 163-177.
\smallskip

\bibitem{Whitt 1993}\textsc{Whitt, W.} (1993). Tail probabilities
with statistical multiplexing and effective bandwidth in
multiclass queues. \emph{Telecommunication Systems}, 2: 71-107.

\smallskip
\bibitem{Whitt 2004}\textsc{Whitt, W.} (2004).
Heavy-traffic limits for loss proportions in single-server queues.
\emph{Queueing Systems}, 46: 507-536.

\smallskip
\bibitem{Whitt 2005} \textsc{Whitt, W.} (2005). Heavy-traffic limits
for the $G/H_2^*/n/m$ queue. \emph{Mathematics of Operations
Research}, 30: 1-27.

\smallskip
\bibitem{Wischik 1999}\textsc{Wischik, D.} (1999). The output of
switch, or, effective bandwidth for networks. \emph{Queueing
Systems}, 32: 383-396.

\smallskip
\bibitem{Wischik 2001}\textsc{Wischik, D.} (2001). Sample path large
deviations for queues with many outputs. \emph{The Annals of
Applied Probability}, 11: 379-404.




\end{thebibliography}
\end{document}